\documentclass[12pt]{article}
\usepackage{amsmath}
\usepackage{amssymb}
\usepackage{amsfonts}
\usepackage{amsthm}
\usepackage{mathrsfs}
\usepackage{units}
\usepackage[all]{xy}

\numberwithin{equation}{section}

\newtheorem{theorem}{Theorem}[section]
\newtheorem{proposition}[theorem]{Proposition}
\newtheorem{corollary}[theorem]{Corollary}
\newtheorem{lemma}[theorem]{Lemma}
\newtheorem{remark}[theorem]{Remark}
\newtheorem{definition}[theorem]{Definition}

\def\F{\mathbb{ F}}
\def\Z{\mathbb{ Z}}

\def\C{\mathbb{ C}}

\def\U{U_q(\mathfrak{g})}
\def\Us{U_q(\widehat{\mathfrak{sl}}_2)}

\date{\small\it June 23, 2006}
\title{Irreducible Modules for the Quantum Affine Algebra $\U$ and its Borel Subalgebra $\U^{\ge 0}$}

\author{ John Bowman\footnote{Supported in part by National Science Foundation Grant \#DMS-0353038.}
}

\begin{document}
\maketitle
\begin{abstract}
We prove a bijection between finite-dimensional irreducible modules for an arbitrary quantum affine algebra $\U$ and finite-dimensional irreducible modules for its Borel subalgebra $\U^{\ge 0}$.
\end{abstract}

\section{Introduction}

Quantum affine algebras have played a significant role in diverse areas of mathematics and physics. These algebras can be associated with vertex operator algebras and feature prominently in conformal field theory (\cite{etingof}, \cite{frankel}). They also have important connections with symmetric functions, in particular with Kostka-Foulkes polynomials (\cite{ck}, \cite{feigin}, \cite{nakaya}). In this paper, we focus on finite-dimensional irreducible modules of an arbitrary quantum affine algebra $\U$. These modules were classified by Chari and Pressley in terms of Drinfeld polynomials in \cite{chari2}. We relate finite-dimensional irreducible modules for $\U$ to finite-dimensional irreducible modules for its Borel subalgebra $\U^{\ge 0}$, where the Borel subalgebra is the nonnegative part of $\U$ with respect to the standard triangular decomposition.

Benkart and Terwilliger showed in \cite{bt} that for the quantum affine algebra $\Us$, there is a bijection between finite-dimensional irreducible modules for $\Us$ and finite-dimensional irreducible modules for the Borel subalgebra $\Us^{\ge 0}$. In particular:
\begin{enumerate}
\item Let $V$ be a finite-dimensional irreducible 
$\Us^{\ge 0}$-module of type $(\varepsilon_1, \varepsilon_2)$. Then the action of $\Us^{\ge 0}$ extends uniquely to an action of $\Us$ on $V$. The resulting 
$\Us$-module structure on $V$ is irreducible and of type $(\varepsilon_1, \varepsilon_2)$.

\item Let $V$ be a finite-dimensional irreducible $\Us$-module of type $(\varepsilon_1, \varepsilon_2)$. When the $\Us$-action is restricted to $\Us^{\ge 0}$, the resulting $\Us^{\ge 0}$-module structure on $V$ is irreducible and of type $(\varepsilon_1, \varepsilon_2)$.
\end{enumerate}
Our main result, contained in Theorem \ref{thm:borel}, extends the result of \cite{bt} to an arbitrary quantum affine algebra $\U$. This situation contrasts with the case where $\mathfrak{g}$ is a finite-dimensional simple Lie algebra. In that case, one can show that every finite-dimensional irreducible $\U^{\ge 0}$-module is one-dimensional.

Our proof contains two key ingredients: the equitable presentation of $\U$ and the notion of a split decomposition. The equitable presentation of $\U$, introduced by Terwilliger in \cite{terwilliger}, has the attractive feature that all of its generators act semisimply on finite-dimensional irreducible $\U$-modules. A split decomposition of a module $V$ is, roughly speaking, a decomposition constructed from two eigenspace decompositions of $V$ by intersecting sums of eigenspaces. For other examples of split decompositions in the literature, we refer the reader to \cite{bt}, \cite{ito}, and \cite{ter2}.

\section{The quantum affine algebra $\U$}

Throughout the paper, we let $\F$ be an algebraically closed field of characteristic 0, and we let $\F^\times$ denote the set of nonzero elements of $\F$. We fix $q \in \F^\times$ such that $q$ is not a root of unity. We let $\Z$ denote the ring of integers, and we set $\nicefrac{1}{2}\, \Z = \{n/2 | n \in \Z\}$.

\begin{definition} \rm
\label{def:A}
Let $n$ denote a positive integer and let $A$ denote a symmetrizable generalized Cartan matrix of order $n$ and affine type \cite[p. 1]{kac}. Since $A$ is symmetrizable there exist relatively prime positive integers $s_1, \ldots, s_n$ such that $s_i A_{ij} = s_j A_{ji}$ for $1 \le i, j \le n$. Since $A$ is of affine type, by \cite[Thm. 4.3]{kac} there exists a column vector $u = (u_1, \ldots, u_n)^T$ such that $u_1, \ldots, u_n$ are relatively prime positive integers and $Au = 0$.

For $1 \le i \le n$ we do the following. We define $q_i = q^{s_i}$. Also, for an integer $m$ we define
\[
[m]_i = \frac{q_i^m - q_i^{-m}}{q_i - q_i^{-1}}
\]
and for $m \ge 0$ we define
\[
[m]_i^{!} = [m]_i [m-1]_i \cdots [2]_i [1]_i.
\]
We interpret $[0]^{!} = 1$. For integers $m \ge r \ge 0$ we define
\[
\genfrac[]{0 cm}{0}mr_i 
= \frac{[m]_i^!}{[r]_i^! [m-r]_i^!}.
\]
Let $\mathfrak{g} = \mathfrak{g}'(A)$ denote the Kac-Moody algebra over $\C$ that corresponds to A \cite[p. xi]{kac}.
\end{definition}

\begin{definition} \rm \cite[p. 280]{chari}
\label{def:U}
With reference to Definition 1.1, $\U$ is the unital associative $\F$-algebra with generators $E_i$, $F_i$, $K_i^{\pm 1}$ $(1 \le i \le n)$ which satisfy the following relations:
\begin{enumerate}
\item[(R1)] $K_i K_i^{-1} = K_i^{-1} K_i = 1$,

\item[(R2)] $K_i K_j = K_j K_i$,

\item[(R3)] $K_i E_j K_i^{-1} = q_i^{A_{ij}} E_j$,

\item[(R4)] $K_i F_j K_i^{-1} = q_i^{-A_{ij}} F_j$,

\item[(R5)] $E_i F_j - F_j E_i = \delta_{ij} \dfrac{K_i - K_i^{-1}}{q_i - q_i^{-1}}$,

\item[(R6)] $\displaystyle\sum_{r=0}^{1 - A_{ij}} (-1)^r
         \genfrac[]{0 cm}{0}{1 - A_{ij}}{r}_i 
         E_i^{1 - A_{ij} - r} E_j E_i^r = 0 \qquad \mbox{if } i \neq j$,

\item[(R7)] $\displaystyle\sum_{r=0}^{1 - A_{ij}} (-1)^r
         \genfrac[]{0 cm}{0}{1 - A_{ij}}{r}_i 
         F_i^{1 - A_{ij} - r} F_j F_i^r = 0 \qquad \mbox{if } i \neq j.$
\end{enumerate}
The expression $\delta_{ij}$ in (R5) is the Kronecker delta.
\end{definition}

We denote by $\U^{\ge 0}$ the subalgebra of $\U$ generated by the elements $E_i, K_i^{\pm 1}$ $(1 \le i \le n)$. We call $\U^{\ge 0}$ the \emph{Borel subalgebra} of $\U$ because of its similarity to the standard Borel subalgebra of the universal enveloping algebra of a finite-dimensional simple Lie algebra over the complex numbers. We want to describe the relationship between finite-dimensional irreducible modules for $\U$ and finite-dimensional irreducible modules for $\U^{\ge 0}$. To this end, we consider finite-dimensional irreducible modules of the following algebra.

\begin{definition} \rm
\label{def:Uge0}
The algebra $U^{\ge 0}$ is the unital associative $\F$-algebra with generators $e_i$, $k_i^{\pm 1}$ $(1 \le i \le n)$ which satisfy the following relations:
\begin{enumerate}
\item[(r1)] $k_i k_i^{-1} = k_i^{-1} k_i = 1$,

\item[(r2)] $k_i k_j = k_j k_i$,

\item[(r3)] $k_i e_j k_i^{-1} = q_i^{A_{ij}} e_j$,

\item[(r4)] $\displaystyle\sum_{r=0}^{1 - A_{ij}} (-1)^r
         \genfrac[]{0 cm}{0}{1 - A_{ij}}{r}_i 
         e_i^{1 - A_{ij} - r} e_j e_i^r = 0 \qquad \mbox{if } i \neq j.$
\end{enumerate}
We let $U^0$ be the subalgebra of $U^{\ge 0}$ generated by $k_i$ $(1 \le i \le n)$. We let $U^{>0}$ be the subalgebra of $U^{\ge 0}$ generated by $e_i$ $(1 \le i \le n)$. Note that by (r3) we have $U^{\ge 0} = U^0 U^{>0} = U^{>0} U^0$.
\label{def:U^0}
\end{definition}

Our first goal is to explain the exact relationship between finite-dimensional irreducible $\U$-modules and finite-dimensional irreducible $U^{\ge 0}$-modules. In order to state our results precisely, it is necessary to make a few comments.

Let $V$ denote a finite-dimensional irreducible $\U$-module. Then for $1 \le i \le n$ the actions of $K_i$ on $V$ are semisimple \cite[Prop. 5.1]{jantzen}. Furthermore (by \cite[Thm. 2.6]{jantzen}), there exist scalars $\varepsilon_i \in \{1, -1\}$ such that each eigenvalue of $K_i$ on $V$ is of the form $\varepsilon_i q_i^m$ $(m \in \Z)$. We call the sequence $\varepsilon = (\varepsilon_i)_{i=1}^n$ the \emph{type} of $V$.

Let $V$ denote a finite-dimensional irreducible $U^{\ge 0}$-module. As we will see in Section 4, for $1 \le i \le n$ the action of $k_i$ on $V$ is semisimple. Moreover, there exist scalars $\alpha_i \in \F^\times$ such that the set of distinct eigenvalues of $k_i$ on $V$ is $\alpha_i \Delta_i$ where $\Delta_i \subseteq \{q_i^m | m \in \Z\}$ and $\Delta_i = \{\theta^{-1} | \theta \in \Delta_i \}$. We refer to the sequence $\alpha = (\alpha_i)_{i=1}^n$ as the \emph{type} of $V$.

Our main results concerning $\U$ and $U^{\ge 0}$ are contained in the following two theorems and the subsequent remark.

\begin{theorem}
\label{thm:Uge0}
Let $V$ denote a finite-dimensional irreducible $U^{\ge 0}$-module of type $\alpha$. Let $\varepsilon \in \{1, -1\}^n$. Then there exists a unique $\U$-module structure on $V$ such that the operators $E_i - \varepsilon_i \alpha_i^{-1} e_i$ and $K_i^{\pm 1} - \varepsilon_i \alpha_i^{\mp 1}k_i^{\pm 1}$ vanish on $V$ for $1 \le i \le n$. This $\U$-module structure is irreducible and of type $\varepsilon$.
\end{theorem}

\begin{theorem}
\label{thm:U}
Let $V$ be a finite-dimensional irreducible $\U$-module of type $\varepsilon$.  Let $\alpha \in (\F^\times)^n$. Then there exists a unique $U^{\ge 0}$-module structure on $V$ such that the operators $E_i - \varepsilon_i \alpha_i^{-1} e_i$ and $K_i^{\pm 1} - \varepsilon_i \alpha_i^{\mp 1}k_i^{\pm 1}$ vanish on $V$ for $1 \le i \le n$. This $U^{\ge 0}$-module structure is irreducible and of type $\alpha$.
\end{theorem}

\begin{remark} \rm
Take $\alpha \in (\F^\times)^n$ and $\varepsilon \in \{1, -1\}^n$. Combining Theorems 1.4 and 1.5, we obtain a bijection between the following two sets:
\begin{enumerate}
\item The isomorphism classes of finite-dimensional irreducible $U^{\ge 0}$-modules of type $\alpha$;

\item The isomorphism classes of finite-dimensional irreducible $\U$-modules of type $\varepsilon$.
\end{enumerate}
\end{remark}

\section{Preliminaries}

For the discussion that follows, we will need the \emph{equitable presentation} of $\U$ introduced by Terwilliger in \cite{terwilliger}.

\begin{proposition} \cite[Thm. 5.1]{terwilliger}
\label{prop:eqU}
The quantum affine algebra $\U$ from Definition \ref{def:U} is isomorphic to the unital associative $\F$-algebra with generators $K_i^{\pm 1}, Y_i, Z_i$ $(1 \le i \le n)$ which satisfy the following relations:
\begin{enumerate}
\item[(E1)] $K_i K_i^{-1} = K_i^{-1} K_i = 1$,

\item[(E2)] $K_i K_j = K_j K_i$,

\item[(E3)] $Y_i K_j - q_i^{-A_{ij}} K_j Y_i = K_i K_j (1 - q_i^{-A_{ij}})$,

\item[(E4)] $Z_i K_j - q_i^{A_{ij}} K_j Z_i = K_i K_j (1 - q_i^{A_{ij}})$,

\item[(E5)] $Z_i Y_i - q_i^2 Y_i Z_i = 1 - q_i^2$,

\item[(E6)] $Z_i Y_j - q_i^{A_{ij}} Y_j Z_i = K_i K_j (1 - q_i^{A_{ij}})
         \quad \mbox{if } i \neq j$,
         
\item[(E7)] $\displaystyle\sum_{r=0}^{1 - A_{ij}} (-1)^r
         \genfrac[]{0 cm}{0}{1 - A_{ij}}{r}_i 
         Y_i^{1 - A_{ij} - r} Y_j Y_i^r =
         K_i^{1 - A_{ij}} K_j \prod_{s=0}^{-A_{ij}} (1 - q_i^{A_{ij} + 2s})
         \quad \mbox{if } i \neq j$,
\item[(E8)] $\displaystyle\sum_{r=0}^{1 - A_{ij}} (-1)^r
         \genfrac[]{0 cm}{0}{1 - A_{ij}}{r}_i 
         Z_i^{1 - A_{ij} - r} Z_j Z_i^r =
         K_i^{1- A_{ij}} K_j \prod_{s=0}^{-A_{ij}} (1 - q_i^{A_{ij} + 2s})
         \quad \mbox{if } i \neq j.$
\end{enumerate}
\begin{tabbing}
12345 \= \kill
The isomorphism with the presentation in Definition \ref{def:U} is \\
\> $K_i^{\pm 1} \rightarrow K_i^{\pm 1}$, \\
\> $Y_i \rightarrow K_i + E_i (q_i-q_i^{-1})$, \\
\> $Z_i \rightarrow K_i - K_i F_i q_i(q_i-q_i^{-1})$. \\
The inverse of this isomorphism is \\
\> $E_i \rightarrow (Y_i - K_i)(q_i-q_i^{-1})^{-1}$, \\
\> $F_i \rightarrow (1 - K_i^{-1} Z_i) q_i^{-1} (q_i - q_i^{-1})^{-1}$, \\
\> $K_i^{\pm 1} \rightarrow K_i^{\pm 1}$.
\end{tabbing}
\end{proposition}

\begin{remark} \rm
The relations in Proposition \ref{prop:eqU} differ slightly from \cite[Thm. 5.1]{terwilliger}. To get our relations from the relations in \cite[Thm. 5.1]{terwilliger}, simply apply the Chevalley involution of $\U$ which sends $K_i^{\pm 1} \rightarrow K_i^{\mp 1}$, $E_i \rightarrow F_i$, and $F_i \rightarrow E_i$.
\end{remark}

For the rest of the paper we identify the two copies of $\U$ given in Definition \ref{def:U} and Proposition \ref{prop:eqU} via the isomorphism in Proposition \ref{prop:eqU}. We now give the equitable presentation of $U^{\ge 0}$.

\begin{proposition} \cite[Thm. 5.1]{terwilliger}
\label{prop:eqUge0}
The algebra $U^{\ge 0}$ from Definition \ref{def:Uge0} is isomorphic to the unital associative $\F$-algebra with generators $k_i^{\pm 1}, y_i$ $(1 \le i \le n)$ which satisfy the following relations:
\begin{enumerate}
\item[(e1)] $k_i k_i^{-1} = k_i^{-1} k_i = 1$,

\item[(e2)] $k_i k_j = k_j k_i$,

\item[(e3)] $y_i k_j - q_i^{-A_{ij}} k_j y_i = k_i k_j (1 - q_i^{-A_{ij}})$,

\item[(e4)] $\displaystyle\sum_{r=0}^{1 - A_{ij}} (-1)^r
        \genfrac[]{0 cm}{0}{1 - A_{ij}}{r}_i 
         y_i^{1 - A_{ij} - r} y_j y_i^r =
        k_i^{1 - A_{ij}} k_j \prod_{s=0}^{-A_{ij}} (1 - q_i^{A_{ij} + 2s})
         \qquad \mbox{if } i \neq j.$
\end{enumerate}
\begin{tabbing}
12345 \= \kill
The isomorphism with the presentation for $U^{\ge 0}$ given in Definition \ref{def:Uge0} is \\
\> $k_i^{\pm 1} \rightarrow k_i^{\pm 1}$, \\
\> $y_i \rightarrow k_i + e_i (q_i-q_i^{-1})$. \\
The inverse of this isomorphism is \\
\> $e_i \rightarrow (y_i - k_i)(q_i-q_i^{-1})^{-1}$, \\
\> $k_i^{\pm 1} \rightarrow k_i^{\pm 1}$.
\end{tabbing}
\end{proposition}

For the rest of the paper we identify the two copies of $U^{\ge 0}$ given in Definition \ref{def:Uge0} and Proposition \ref{prop:eqUge0} via the isomorphism in Proposition \ref{prop:eqUge0}. We will make use of the following additional relation in $U^{\ge 0}$.

\begin{lemma} \cite[line (19)]{terwilliger}
\label{lem:2.1}
The following relation holds in $U^{\ge 0}$ for all $1 \le i, j \le n$ with $i \neq j$:
\begin{equation}
\label{eq:2.1}
\displaystyle\sum_{r=0}^{1 - A_{ij}} (-1)^r
         \genfrac[]{0 cm}{0}{1 - A_{ij}}{r}_i 
         y_i^{1 - A_{ij} - r} e_j y_i^r = 0.
\end{equation}
\end{lemma}

\begin{proof}
This is proved for $\U$ in \cite[pp. 309--311]{terwilliger}. The proof carries over to $U^{\ge 0}$ in a straightforward fashion.
\end{proof}

\section{Eigenspace decompositions}

Throughout this section, we fix a finite-dimensional irreducible $U^{\ge 0}$-module $V$. By a {\it decomposition of $V$} we mean a sequence of subspaces of $V$ whose direct sum is $V$.

\begin{lemma}
\label{lem:U_i}
For $1 \le i \le n$ there exist $\alpha_i \in \F^\times$, a positive $d_i \in \nicefrac{1}{2} \, \Z$, and a decomposition $\{U_i(s)\}_{s=0}^{2d_i}$ of $V$ such that the following hold:
\begin{enumerate}
\item[(i)] $(k_i - \alpha_i q_i^{s-d_i}I) U_i(s) = 0$ for $0 \le s \le 2d_i$,

\item[(ii)] $U_i(0) \neq 0$ and $U_i(2d_i) \neq 0$.
\end{enumerate}
Moreover, $k_i$ is semisimple on $V$.
\end{lemma}

\begin{proof}
For $\theta \in \F$ let $U(\theta) = \{v \in V | k_i v = \theta v\}$. From (r2), (r3) we observe the following for $1 \le j \le n$:
\begin{equation}
\label{eq:3.1}
k_j U(\theta) \subseteq U(\theta),
\end{equation}
\begin{equation}
\label{eq:3.2}
e_j U(\theta) \subseteq U(\theta q_i^{A_{ij}}).
\end{equation}
Since $\F$ is algebraically closed, $k_i$ has an eigenvalue in $\F$. Thus there exists $\eta \in \F$ such that $U(\eta) \neq 0$. Since $k_i$ is invertible on $V$, we have $\eta \neq 0$. The scalars $\eta q_i^k$ $(k \in \Z)$ are mutually distinct since $q_i$ is not a root of unity. Because $V$ is finite-dimensional, there exists $m \in \Z$ such that $U(\eta q_i^{m}) \neq 0$ and $U(\eta q_i^s) = 0$ if $s < m$. Also, there exists $M \in \Z$ such that $U(\eta q_i^{M}) \neq 0$ and $U(\eta q_i^s) = 0$ if $s > M$. Set $\alpha_i := \eta q_i^{(m+M)/2}$ and $d_i := (M - m)/2$. Set $U_i(s) := U(\alpha_i q_i^{s - d_i})$ for $0 \le s \le 2d_i$. By the construction we have (i) and (ii). Moreover, the sum
$\sum_{s=0}^{2d_i} U_i(s)$
is invariant under $U^{\ge 0}$ by (\ref{eq:3.1}), (\ref{eq:3.2}). The sum is nonzero since $U_i(0) \neq 0$, so by irreducibility this sum is equal to $V$. Since the $U_i(s)$ are eigenspaces for $k_i$ corresponding to distinct eigenvalues, the sum is direct. Thus $\{U_i(s)\}_{s=0}^{2d_i}$ is a decomposition of $V$ as desired. Clearly $k_i$ is semisimple on $V$.
\end{proof}

\begin{definition} \rm
For notational convenience, we define $U_i(s) = 0$ if $s < 0$ or $s > 2d_i$ $(1 \le i \le n)$.
\end{definition}

\begin{definition} \rm
With reference to Lemma \ref{lem:U_i}, the sequence $\alpha = (\alpha_i)_{i=1}^n$ is the \emph{type} of $V$. We refer to the sequence $d = (d_i)_{i=1}^n$ as the \emph{shape} of $V$. 
\end{definition}

\begin{remark} \rm
In Section 4 we will show $d_i \in \Z$ for $1 \le i \le n$.
\end{remark}

\begin{lemma}
\label{lem:U_i.act}
For $1 \le i \le n$ let the decomposition $\{U_i(s)\}_{s=0}^{2d_i}$ be as in Lemma \ref{lem:U_i}. Then the following hold for $1 \le j \le n$ and $0 \le s \le 2d_i$:
\begin{enumerate}
\item[(i)] $k_j U_i(s) = U_i(s)$,

\item[(ii)] $e_j U_i(s) \subseteq U_i(s + A_{ij})$.
\end{enumerate}
\end{lemma}

\begin{proof}
(i) By (\ref{eq:3.1}), $k_i U_i(s) \subseteq U_i(s)$. The result follows since $k_i$ is invertible on $V$. (ii) follows immediately from (\ref{eq:3.2}).
\end{proof}

\begin{lemma}
\label{lem:S}
Let $S$ be a subspace of $V$ that is invariant under $U^{>0}$. 
Suppose there exist integers $i$ and $s$ such that $0 \neq U_i(s) \subseteq S$. Then $S = V$.
\end{lemma}
\begin{proof}
By hypothesis we have $U^{> 0} U_i(s) \subseteq U^{> 0} S \subseteq S$. Thus it suffices to show $U^{> 0} U_i(s) = V$. By construction, $U^{>0} U_i(s)$ is nonzero and invariant under $U^{>0}$. Using (r3) and Lemma \ref{lem:U_i.act}(i) we find $U^{> 0} U_i(s)$ is invariant under $U^0$. Thus $U^{>0} U_i(s)$ is invariant under $U^{\ge 0}$, and we have $U^{>0}U_i(s) = V$ since $V$ is an irreducible $U^{\ge 0}$-module. The result follows.
\end{proof}

\begin{lemma}
For $1 \le i \le n$ let the decomposition $\{U_i(s)\}_{s=0}^{2d_i}$ be as in Lemma \ref{lem:U_i}. Then the following holds for $0 \le  s \le 2d_i$:
\begin{equation}
\label{eq:y_i.U_i}
(y_i - \alpha_i q_i^{s - d_i}I) U_i(s) \subseteq U_i(s+2).
\end{equation}
\end{lemma}

\begin{proof}
In (\ref{eq:y_i.U_i}) evaluate the left-hand side using $y_i = k_i + (q_i - q_i^{-1})e_i$. The result follows from Lemma \ref{lem:U_i}(i) and Lemma \ref{lem:U_i.act}(ii).
\end{proof}

\begin{lemma}
\label{lem:V_i}
For $1 \le i \le n$ there exists a decomposition $\{V_i(s)\}_{s=0}^{2d_i}$ of $V$ such that
\begin{equation}
\label{eq:V_i}
(y_i - \alpha_i q_i^{s - d_i}I) V_i(s) = 0 \qquad (0 \le s \le 2d_i).
\end{equation}
Moreover, $y_i$ is semisimple on $V$.
\end{lemma}

\begin{proof} By (\ref{eq:y_i.U_i}) the product $\prod_{s=0}^{2d_i} (y_i - \alpha_i q_i^{s - d_i} I)$ vanishes on $V$. Since the scalars $\alpha_i q_i^{s - d_i}$ are mutually distinct, we find that $y_i$ is semisimple on $V$ with eigenvalues contained in the set $\{ \alpha_i q_i^{s - d_i} | 0 \le s \le 2d_i \}$. We define $V_i(s) := \{v \in V |  y_i v = \alpha_i q_i^{s - d_i}\}$.
Thus (\ref{eq:V_i}) holds by construction. Since $y_i$ is semisimple with eigenspaces $V_i(s)$, we see that $\{V_i(s)\}_{s=0}^{2d_i}$ is a decomposition of $V$ as desired.
\end{proof}

\begin{definition} \rm
For notational convenience, we define $V_i(s) = 0$ if $s < 0$ or $s > 2d_i$ $(1 \le i \le n)$.
\end{definition}

\begin{lemma}
\label{lem:V_i.act}
For $1 \le i \le n$ let the decomposition $\{V_i(s)\}_{s=0}^{2d_i}$ be as in Lemma \ref{lem:V_i}. Then the following hold for $0 \le s \le 2d_i$:
\begin{enumerate}
\item[(i)] $(k_i^{-1} - \alpha_i^{-1} q_i^{d_i - s}I) V_i(s) \subseteq V_i(s+2)$,

\item[(ii)] $e_j V_i(s) \subseteq \displaystyle\sum_{t=0}^{-A_{ij}} V_i(s + A_{ij} + 2t)$ if $i \neq j$ \quad $(1 \le j \le n)$.
\end{enumerate}
\end{lemma}

\begin{proof} (i) By (e3) (with $i=j$) we have 
\[
(y_i k_i^{-1} - q_i^2 k_i^{-1} y_i - I + q_i^2I) V_i(s) = 0.
\]
Applying (\ref{eq:V_i}) we obtain
\begin{align*}
0 &= (y_i k_i^{-1} - \alpha_i q_i^{s - d_i + 2} k_i^{-1} - \alpha_i^{-1} q_i^{d_i - s} y_i + q_i^2I) V_i(s) \\
&= (y_i - \alpha_i q_i^{s - d_i + 2}) (k_i^{-1} - \alpha_i^{-1} q_i^{d_i - s}I) V_i(s).
\end{align*}
The result follows.

(ii) Let $S$ denote the left-hand side of (\ref{eq:2.1}) and note that $S=0$. Applying $S$ to $V_i(s)$ and using (\ref{eq:V_i}), we find that $S$ agrees with
\begin{equation}
\label{eq:V_i.act.1}
\prod_{t=0}^{-A_{ij}} (y_i - \alpha_i q_i^{s-d_i+A_{ij}+2t}I)e_j
\end{equation}
on $V_i(s)$. Therefore (\ref{eq:V_i.act.1}) is zero on $V_i(s)$. The result follows using (\ref{eq:V_i}).
\end{proof}

\begin{lemma}
\label{lem:UVsum}
For $1 \le i \le n$ let the decompositions $\{U_i(s)\}_{s=0}^{2d_i}$ and $\{V_i(s)\}_{s=0}^{2d_i}$ be as in Lemmas \ref{lem:U_i} and \ref{lem:V_i} respectively. Then for $0 \le s \le 2d_i$ we have
\[
\sum_{t = 0}^{\infty} U_i(s + 2t) = \sum_{t=0}^{\infty} V_i(s + 2t).
\]
\end{lemma}

\begin{proof} We set $X := \sum_{t=0}^{\infty} U_i(s + 2t)$ and $X' := \sum_{t=0}^{\infty} V_i(s + 2t)$. We show that $X = X'$.
Let $N$ be an integer with $N > d_i - s/2$.
Note that $X = \sum_{t=0}^N U_i(s + 2t)$ and $X' = \sum_{t=0}^N V_i(s + 2t)$.
Let $T:= \prod_{t=0}^N (y_i - \alpha_i q_i^{s - d_i + 2t}I)$.
Then $X' = \{v \in V | T v = 0\}$ by (\ref{eq:V_i}), and $T X = 0$ by (\ref{eq:y_i.U_i}). Thus $X \subseteq X'$.
Now let
\[
S := \prod_{r=0}^{s-1} (y_i - \alpha_i q_i^{r - d_i}I)
\prod_{t=0}^N (y_i - \alpha_i q_i^{s - d_i + 2t + 1}I).
\]
Observe that $S V = X'$ by (\ref{eq:V_i}), and $S V \subseteq X$ by (\ref{eq:y_i.U_i}), so $X' \subseteq X$. Thus $X' = X$, completing the proof. \end{proof}

\section{Split decompositions}

We fix a finite-dimensional irreducible $U^{\ge 0}$-module $V$ of type $\alpha$ and shape $d$. For $1 \le i \le n$ we use the decomposition $\{U_i(s)\}_{s=0}^{2d_i}$ of $V$ from Lemma \ref{lem:U_i} and the decomposition $\{V_i(s)\}_{s=0}^{2d_i}$ of $V$ from Lemma \ref{lem:V_i} to construct a third decomposition $\{W_i(s)\}_{s=0}^{2d_i}$ of $V$.

\begin{definition} \rm Let $s$ be an integer. For $1 \le i \le n$ we define
\label{def:W_i}
\[
W_i(s) = \left( \sum_{t=0}^{\infty} U_i(s - 2t) \right) \cap 
\left( \sum_{t=0}^{\infty} V_i(2d_i - s - 2t) \right).
\]
\end{definition}

\noindent Note that $W_i(s) = 0$ if $s < 0$ or $s > 2d_i$. We will show that $\{W_i(s)\}_{s=0}^{2d_i}$ is a decomposition of $V$. Toward this end, the following definition will be useful.

\begin{definition} \rm
\label{def:W_i2}
Let $u$ and $v$ be integers. For $1 \le i \le n$ we define
\[
W_i(u,v) = \left( \sum_{t=0}^{\infty} U_i(u - 2t) \right) \cap 
\left( \sum_{t=0}^{\infty} V_i(2d_i - v - 2t) \right).
\]
\end{definition}

\begin{lemma}
\label{lem:W_i2.act}
Choose integers $u$ and $v$ with $0 \le u, v \le 2d_i$. With reference to Definition \ref{def:W_i2}, the following hold for $1 \le i \le n$:
\begin{enumerate}
\item[(i)] $(k_i^{-1} - \alpha_i^{-1} q_i^{d_i - u}I) W_i(u,v) \subseteq W_i(u-2, v-2)$,

\item[(ii)] $(y_i - \alpha_i q_i^{d_i - v}I) W_i(u,v) \subseteq W_i(u+2,v+2)$,

\item[(iii)] $e_j W_i(u,v) \subseteq W_i(u+A_{ij}, v+A_{ij})$ if $i \neq j$.
\end{enumerate}
\end{lemma}

\begin{proof}
(i) By Lemma \ref{lem:U_i}(i) we have
\[
(k_i^{-1} - \alpha_i^{-1} q_i^{d_i - u} I)
\sum_{t=0}^{\infty} U_i(u - 2t)  \subseteq
\sum_{t=0}^{\infty} U_i(u - 2 - 2t).
\]
By Lemma \ref{lem:V_i.act}(i) we have
\[
(k_i^{-1} - \alpha_i^{-1} q_i^{d_i - u} I)
\sum_{t=0}^{\infty} V_i(2d_i - v - 2t)  \subseteq
\sum_{t=0}^{\infty} V_i(2d_i - v + 2 - 2t).
\]
The result follows from the above inclusions and Definition \ref{def:W_i2}.

(ii) By (\ref{eq:y_i.U_i}) we have
\[
(y_i - \alpha_i q_i^{d_i - v}I)
\sum_{t=0}^{\infty} U_i(u - 2t)  \subseteq
\sum_{t=0}^{\infty} U_i(u + 2 - 2t).
\]
By (\ref{eq:V_i}) we have
\[
(y_i - \alpha_i q_i^{d_i - v} I)
\sum_{t=0}^{\infty} V_i(2d_i - v - 2t)  \subseteq
\sum_{t=0}^{\infty} V_i(2d_i - v - 2 - 2t).
\]
The result follows from the above inclusions and Definition \ref{def:W_i2}.

(iii) By Lemma \ref{lem:U_i.act}(ii) we have
\[
e_j \sum_{t=0}^{\infty} U_i(u - 2t) \subseteq
\sum_{t=0}^{\infty} U_i(u + A_{ij} - 2t).
\]
By Lemma \ref{lem:V_i.act}(ii) we have
\[
e_j \sum_{t=0}^{\infty} V_i(2d_i - v - 2t) \subseteq
\sum_{t=0}^{\infty} V_i(2d_i - v - A_{ij} - 2t).
\]
The result follows from the above inclusions and Definition \ref{def:W_i2}.
\end{proof}

\begin{corollary}
\label{cor:Winv}
With reference to Definition \ref{def:W_i2} and Definition \ref{def:U^0}, the sum \\
$\sum_{s \in \Z} W_i(u+s,v+s)$ is invariant under $U^{> 0}$.
\end{corollary}

\begin{proof}
Set $S := \sum_{s \in \Z} W_i(u+s,v+s)$. For $1 \le j \le n$ with $j \neq i$ we have $e_j S \subseteq
S$ by Lemma \ref{lem:W_i2.act}(iii). By Lemma \ref{lem:W_i2.act}(i) we have $k_i^{-1} S \subseteq S$, so $k_i S \subseteq S$ also. By Lemma \ref{lem:W_i2.act}(ii) we have $y_i S \subseteq S$. Using $e_i = (y_i - k_i)(q_i - q_i^{-1})^{-1}$ we have $e_i S \subseteq S$. The result follows.
\end{proof}

\begin{lemma}
\label{lem:u<v}
Let $u$ and $v$ be integers with $u < v$. Then for $1 \le i \le n$ the space $W_i(u,v)$ from Definition \ref{def:W_i2} is zero.
\end{lemma}

\begin{proof}
We assume $W_i(u,v) \neq 0$ and derive a contradiction. Set $S := \sum_{s \in \Z} W_i(u+s,v+s)$ and note $S \neq 0$.
We first show $U^0 S = V$. Clearly $U^0 S$ is nonzero and invariant under $U^0$. By (r3) and Corollary \ref{cor:Winv}, $U^0 S$ is invariant under $U^{>0}$. Thus $U^0 S$ is a nonzero subspace invariant under $U^{\ge 0}$, and we have $U^0 S = V$ since $V$ is an irreducible $U^{\ge 0}$-module.
Now we show $U^0 S \neq V$ for a contradiction.
To this end, we claim that $U^0 S$ is contained in the sum
\begin{equation}
\label{eq:u<v}
\sum_{s=0}^{2d_i + u - v} U_i(s).
\end{equation}
For $s \in \Z$ note that $W_i(u+s,v+s)$ is contained in (\ref{eq:u<v}) by Definition \ref{def:W_i2}. By this and since each $U_i(s)$ is $U^0$-invariant by Lemma \ref{lem:U_i.act}, we find that $U^0 S$ is contained in (\ref{eq:u<v}). But (\ref{eq:u<v}) is not equal to $V$ by Lemma \ref{lem:U_i}. Thus $U^0 S \neq V$, and we have a contradiction. The result follows.
\end{proof}

\begin{theorem}
\label{thm:W_i}
Choose $i$ with $1 \le i \le n$. With the notation of Definition 4.1, $\{W_i(s)\}_{s=0}^{2d_i}$ is a decomposition of $V$.
Further, $d_i \in \Z$.
\end{theorem}

\begin{proof}
We first show $d_i \in \Z$. Suppose not. By construction, $d_i \in \nicefrac{1}{2}\,\Z$, so $2d_i$ is an odd integer. Now applying Definition \ref{def:W_i2} we find
\begin{align*}
W_i(0, 1) &= U_i(0) \cap \bigl( V_i(2d_i - 1) + \cdots + V_i(2) + V_i(0) \bigr) \\
&= U_i(0) \cap \bigl( U_i(2d_i - 1) + \cdots + U_i(2) + U_i(0) \bigr) \mbox{ (by Lemma \ref{lem:UVsum})} \\
&= U_i(0).
\end{align*}
Observe $W_i(0, 1) = 0$ by Lemma \ref{lem:u<v} and $U_i(0) \neq 0$ by Lemma \ref{lem:U_i}(ii) for a contradiction. Thus $d_i \in \Z$.

We now show $\{W_i(s)\}_{s=0}^{2d_i}$ is a decomposition of $V$. To this end, consider the sum $S := \sum_{s=0}^{2d_i} W_i(s)$. We first show $S = V$. Note $W_i(s) = W_i(s,s)$ by Definition \ref{def:W_i2}. Thus $S$ is invariant under $U^{>0}$ by Corollary \ref{cor:Winv}. 
We have
\begin{align*}
W_i(0) &= U_i(0) \cap \bigl( V_i(2d_i) + \cdots + V_i(2) + V_i(0) \bigr) \\
&= U_i(0) \cap \bigl( U_i(2d_i) + \cdots + U_i(2) + U_i(0) \bigr) \mbox{ (by Lemma \ref{lem:UVsum})} \\
&= U_i(0)
\end{align*}
so $U_i(0) = W_i(0) \subseteq S$. Thus $S=V$ by Lemma \ref{lem:S}.

Finally, we show that the sum $\sum_{s=0}^{2d_i} W_i(s)$ is direct. Since this sum is equal to $V$, $V$ is the sum of
\begin{equation}
\label{eq:W_i.1}
\sum_{s=0}^{d_i} W_i(2s)
\end{equation}
and
\begin{equation}
\label{eq:W_i.2}
\sum_{s=0}^{d_i-1} W_i(2s+1).
\end{equation}
By Definition \ref{def:W_i}, (\ref{eq:W_i.1}) is contained in
$\sum_{s=0}^{d_i} U_i(2s)$ and (\ref{eq:W_i.2}) is contained in \\
$\sum_{s=0}^{d_i-1} U_i(2s+1)$. Since $\{U_i(s)\}_{s=0}^{2d_i}$ is a decomposition of $V$, the intersection of (\ref{eq:W_i.1}) and (\ref{eq:W_i.2}) is zero. Thus $V$ is the direct sum of (\ref{eq:W_i.1}) and (\ref{eq:W_i.2}). We now show that the sum (\ref{eq:W_i.1}) is direct. It suffices to show that
\[
\bigl( W_i(0) + W_i(2) + \cdots + W_i(2m) \bigr) \cap W_i(2m+2) = 0
\]
for $0 \le m < d_i$.
By Definition \ref{def:W_i} we have
\begin{align*}
\bigl( W_i(0)& + W_i(2) + \cdots + W_i(2m) \bigr) \cap W_i(2m+2) \\
&\subseteq \left( \sum_{t=0}^{\infty} U_i(2m -2t) \right) \cap \left( \sum_{t=0}^{\infty} V_i(2d_i -2m - 2 -2t) \right) \\
&= W_i(2m,2m+2) = 0 \mbox{ (by Lemma \ref{lem:u<v})},
\end{align*}
so the sum (\ref{eq:W_i.1}) is direct. Similarly we see the sum (\ref{eq:W_i.2}) is direct. The result follows.
\end{proof}

\begin{lemma}
\label{lem:W_i.act}
For $1 \le i \le n$ let the decomposition $\{W_i(s)\}_{s=0}^{2d_i}$ be as in Definition \ref{def:W_i}. Then the following hold for $0 \le s \le 2d_i$:
\begin{enumerate}
\item[(i)] $(k_i^{-1} - \alpha_i^{-1} q_i^{d_i - s}I) W_i(s) \subseteq W_i(s-2)$,

\item[(ii)] $(y_i - \alpha_i q_i^{d_i - s}I) W_i(s) \subseteq W_i(s+2)$,

\item[(iii)] $e_j W_i(s) \subseteq W_i(s+A_{ij})$ if $i \neq j$.
\end{enumerate}
\end{lemma}

\begin{proof}
Note that $W_i(s) = W_i(s,s)$, and apply Lemma \ref{lem:W_i2.act}.
\end{proof}

\begin{lemma}
\label{lem:Wsum}
For $1 \le i \le n$ let the decompositions $\{U_i(s)\}_{s=0}^{2d_i}$ and $\{V_i(s)\}_{s=0}^{2d_i}$ be as in Lemmas \ref{lem:U_i} and \ref{lem:V_i} respectively. Let the decomposition $\{W_i(s)\}_{s=0}^{2d_i}$ be as in Definition \ref{def:W_i}. Then the following hold for $0 \le s \le 2d_i$:
\begin{enumerate}
\item[(i)] $\displaystyle\sum_{t=0}^{\infty} W_i(s - 2t) = \displaystyle\sum_{t=0}^{\infty} U_i(s - 2t)$,

\item[(ii)] $\displaystyle\sum_{t=0}^{\infty} W_i(s + 2t) = \displaystyle\sum_{t=0}^{\infty} V_i(2d_i - s - 2t)$.
\end{enumerate}
\end{lemma}

\begin{proof}
(i)  Let $X := \sum_{t=0}^{\infty} W_i(s - 2t)$ and $X' := \sum_{t=0}^{\infty} U_i(s - 2t)$.
Note that $X \subseteq X'$ by Definition \ref{def:W_i}. To see the reverse inclusion, choose an integer $N$ with $N > s/2$. Define
\[
S := \prod_{r=-d_i}^{d_i-s-1} (k_i^{-1} - \alpha_i^{-1} q_i^{r}I)
\prod_{t=0}^{N} (k_i^{-1} - \alpha_i^{-1} q_i^{d_i - s + 2t+1}I).
\]
By Lemma \ref{lem:U_i}(i) we have $SV = X'$. By Lemma \ref{lem:W_i.act}(i) we have $S V \subseteq X$. Thus $X' \subseteq X$, and $X = X'$ as desired.

(ii) Let $Y := \sum_{t=0}^{\infty} W_i(s + 2t)$, and $Y' := \sum_{t=0}^{\infty} V_i(2d_i - s - 2t)$. Note that $Y \subseteq Y'$ by Definition \ref{def:W_i}.  To see the reverse inclusion, choose an integer $N$ with $N > d_i - s/2$. Define
\[
T := \prod_{r=d_i-s+1}^{d_i} (y_i - \alpha_i q_i^{r}I)
\prod_{t=0}^N (y_i - \alpha_i q_i^{d_i - s - 2t -1}I).
\]
By (\ref{eq:V_i}) we have $TV = Y'$. By Lemma \ref{lem:W_i.act}(ii) we have $TV \subseteq Y$. Thus $Y' \subseteq Y$, and $Y = Y'$ as desired.
\end{proof}

\section{The $\U$-module structure}

We fix a finite-dimensional irreducible $U^{\ge 0}$-module $V$ of type $\alpha$ and shape $d$. Let $\varepsilon = (\varepsilon_i)_{i=1}^n \in \{-1, 1\}^n$. In this section, we will construct a $\U$-module structure on $V$ of type $\varepsilon$ and show that this structure satisfies the requirements of Theorem \ref{thm:Uge0}.

\begin{definition} \rm
\label{def:U.act}
Recall the decomposition $\{W_i(s)\}_{s=0}^{2d_i}$ from Definition 4.1. With reference to Proposition \ref{prop:eqU}, we let the generators $K_i^{\pm 1}$, $Y_i$, and $Z_i$ of $\U$ act on $V$ as follows for $1 \le i \le n$.
Let $K_i^{\pm 1}$ act as $\varepsilon_i \alpha_i^{\mp 1} k_i^{\pm 1}$. Let $Y_i$ act as $\varepsilon_i \alpha_i^{-1} y_i$. Let $Z_i$ act as the linear operator $V \rightarrow V$ such that for $0 \le s \le 2d_i$, $W_i(s)$ is the eigenspace of $Z_i$ with eigenvalue $\varepsilon_i q_i^{s - d_i}$.
\end{definition}

We now show that these actions satisfy the relations of $\U$, giving a $\U$-module structure on $V$.

\begin{remark} \rm
\label{rem:E_i.act}
Recall that $E_i = (Y_i - K_i)q_i^{-1} (q_i - q_i^{-1})^{-1}$ and $e_i = (y_i - k_i)q_i^{-1} (q_i - q_i^{-1})^{-1}$. In view of Definition \ref{def:U.act}, $E_i$ acts on $V$ as $\varepsilon_i \alpha_i^{-1} e_i$.
\end{remark}

\begin{lemma}
\label{lem:U.rel}
With reference to Definition \ref{def:U.act} and Remark \ref{rem:E_i.act}, the following relations hold on $V$ for $1 \le i \le n$:
\begin{enumerate}
\item[(i)] $Z_i K_i - q_i^2 K_i Z_i = (1 - q_i^2) K_i^2$,

\item[(ii)] $Z_i Y_i - q_i^2 Y_i Z_i = (1 - q_i^2) I$,

\item[(iii)] $Z_i E_h = q_i^{A_{ih}} E_h Z_i$ if $i \neq h$.
\end{enumerate}
\end{lemma}

\begin{proof}
Since $\{W_i(s)\}_{s=0}^{2d_i}$ is a decomposition of $V$, it suffices to show that (i)--(iii) hold on $W_i(s)$ for $0 \le s \le 2d_i$.

(i) 
By Lemma \ref{lem:W_i.act}(i) we have the following for $0 \le s \le 2d_i$:
\begin{align*}
0 &= (Z_i - \varepsilon_i q_i^{s - d_i - 2}I) (K_i^{-1} - \varepsilon_i q_i^{d_i - s}I) W_i(s) \\
&= (Z_i K_i^{-1} - \varepsilon_i q_i^{s - d_i - 2}K_i^{-1} - \varepsilon_i q_i^{d_i - s} Z_i + q_i^{-2}I) W_i(s) \\
&= \bigl( Z_i K_i^{-1} - q_i^{-2} K_i^{-1} Z_i + (q_i^{-2} - 1)I \bigr) W_i(s).
\end{align*}
The result follows after rearrangement of terms.

(ii) 
By Lemma \ref{lem:W_i.act}(ii) we have the following for $0 \le s \le 2d_i$:
\begin{align*}
0 &= (Z_i - \varepsilon_i q_i^{s - d_i + 2}I) (Y_i - \varepsilon_i q_i^{d_i - s}I) W_i(s) \\
&= (Z_i Y_i - \varepsilon_i q_i^{s-d_i+2}Y_i -\varepsilon_i q_i^{d_i-s} Z_i + q_i^2I)W_i(s)  \\
&= \bigl( Z_i Y_i - q_i^2 K_i Z_i + (q_i^2 - 1)I \bigr) W_i(s).
\label{eq:5.2}
\end{align*}
The result follows.

(iii) 
By Lemma \ref{lem:W_i.act}(iii) we have the following for $0 \le s \le 2d_i$:
\begin{align*}
0 &= (Z_i - \varepsilon_i q_i^{s - d_i + A_{ih}}I) E_h W_i(s) \\
&= (Z_i E_h - q_i^{A_{ih}} E_h Z_i) W_i(s).
\end{align*}
The result follows.
\end{proof}

\begin{lemma}
\label{lem:K_j2}
Choose $i$ and $j$ with $1 \le i,j \le n$. With the notation of Definition \ref{def:U.act}, the following relation holds on $V$:
\begin{equation}
\label{eq:K_j}
Z_i K_j - q_i^{A_{ij}} K_j Z_i = K_i K_j (1 - q_i^{A_{ij}}).
\end{equation}
\end{lemma}

\begin{proof}
The result holds for $i=j$ by Lemma \ref{lem:U.rel}(i), so we may assume $i \neq j$.
It suffices to show $Z_i - q_i^{A_{ij}} K_j Z_i K_j^{-1} = (1 - q_i^{A_{ij}}) K_i$ on $V$. To this end, we define $S := \{ v \in V | (Z_i - q_i^{A_{ij}} K_j Z_i K_j^{-1})v = (1 - q_i^{A_{ij}}) K_i v \}$ and show $S = V$. We do this using Lemma \ref{lem:S} with $s = 0$. Recall $U_i(0) \neq 0$ by Lemma \ref{lem:U_i}(ii). We claim that $U_i(0) \subseteq S$. By Lemma \ref{lem:Wsum}, $U_i(0)$ is equal to $W_i(0)$, the eigenspace of $Z_i$ with eigenvalue $\varepsilon_i q_i^{-d_i}$. By Lemma \ref{lem:U_i.act}(i), $U_i(0)$ is equal to $K_j U_i(0) = K_j W_i(0)$, the eigenspace of $K_j Z_i K_j^{-1}$ with eigenvalue $\varepsilon_i q_i^{-d_i}$. Since $U_i(0)$ is also the eigenspace of $K_i$ with eigenvalue $\varepsilon_i q_i^{-d_i}$, we have $(Z_i - q_i^{A_{ij}} K_j Z_i K_j^{-1})v = (1 - q_i^{A_{ij}})K_i v$ for $v \in U_i(0)$, so $U_i(0) \subseteq S$ as claimed. We now show that $S$ is invariant under $U^{>0}$. To this end, we first show
\begin{equation}
\label{eq:E_h}
(Z_i - q_i^{A_{ij}} K_j Z_i K_j^{-1}) E_h = q_i^{A_{ih}} E_h (Z_i - q_i^{A_{ij}} K_j Z_i K_j^{-1})
\end{equation}
for $1 \le h \le n$. To verify (\ref{eq:E_h}) for $h \neq i$, 
take the term $E_h$ that appears on the left-hand side and routinely pull it to the left using (r3) and Lemma \ref{lem:U.rel}(iii).
For $h = i$,  it suffices to show (\ref{eq:E_h}) with $E_h$ replaced by $Y_i - K_i$ in view of Remark \ref{rem:E_i.act}.
To verify the modified (\ref{eq:E_h}), take the factor $Y_i - K_i$ that appears on the left-hand side and pull it to the left using (e2), (e3), and Lemma \ref{lem:U.rel}(i), (ii).
It follows from (\ref{eq:E_h}), Remark \ref{rem:E_i.act}, and (r3) that $S$ is invariant under $U^{> 0}$. We have now shown that $S$ contains $U_i(0) \neq 0$ and is invariant under $U^{>0}$, so by Lemma \ref{lem:S} we have $S = V$. The result follows.
\end{proof}

\begin{lemma}
\label{lem:Y_j}
Choose $i$ and $j$ with $1 \le i,j \le n$ and $i \neq j$. With the notation of Definition \ref{def:U.act}, the following relation holds on $V$:
\begin{equation}
\label{eq:Y_j}
Z_i Y_j - q_i^{A_{ij}} Y_j Z_i = K_i K_j (1 - q_i^{A_{ij}}).
\end{equation}
\end{lemma}

\begin{proof}
By Lemma \ref{lem:U.rel}(iii) and Remark \ref{rem:E_i.act} we have
\[
Z_i (Y_j - K_j) - q_i^{A_{ij}} (Y_j - K_j) Z_i = 0.
\]
The result follows using (\ref{eq:K_j}).
\end{proof}

\begin{lemma}
\label{lem:Z_j.W_i}
Choose $i$ and $j$ with $1 \le i,j \le n$ and $i \neq j$. Let the decomposition $\{W_i(s)\}_{s=0}^{2d_i}$ be as in Definition \ref{def:W_i}. With reference to Definition \ref{def:U.act}, we have the following for $0 \le s \le 2d_i$:
\begin{equation}
\label{eq:Z_j.W_i}
(Z_j - K_j) W_i(s) \subseteq \sum_{t=0}^{-A_{ij}} W_i(s + A_{ij} + 2t).
\end{equation}
\end{lemma}

\begin{proof}
Exchanging $i$ and $j$ in (\ref{eq:K_j}) and applying Lemma \ref{lem:U_i}(i) we have
\begin{align*}
0 &= (\varepsilon_i q_i^{s - d_i} Z_j - q_i^{A_{ij}} K_i Z_j - \varepsilon_i q_i^{s-d_i} K_j + q_i^{A_{ij}} K_i K_j) U_i(s) \\
&= (K_i - \varepsilon_i q_i^{s - d_i - A_{ij}}) (Z_j - K_j) U_i(s).
\end{align*}
Thus $(Z_j - K_j) U_i(s) \subseteq U_i(s - A_{ij})$. We also have
\begin{align*}
0 &= (Z_j Y_i - q_i^{A_{ij}} Y_i Z_j - K_jY_i + q_i^{A_{ij}}Y_iK_j) V_i(s) \qquad \mbox{(by (\ref{eq:Y_j}) and (e3))} \\
&= (Y_i - \varepsilon_i q_i^{s - d_i - A_{ij}}I)(Z_j - K_j)V_i(s) \qquad \mbox{(by (\ref{eq:V_i}) and Definition \ref{def:U.act})}.
\end{align*}
Thus $(Z_j - K_j) V_i(s) \subseteq V_i(s - A_{ij})$. It follows that
\begin{align*}
(Z_j - K_j) W_i(s) &=
(Z_j - K_j) \left( \sum_{t=0}^{\infty} U_i(s-2t) \right)
\cap \left( \sum_{t=0}^{\infty} V_i(2d_i - s - 2t) \right) \\
& \subseteq \left( \sum_{t=0}^{\infty} U_i(s-2t-A_{ij}) \right)
\cap \left( \sum_{t=0}^{\infty} V_i(2d_i - s - 2t - A_{ij}) \right) \\
& =  \left( \sum_{t=0}^{\infty} W_i(s-2t-A_{ij}) \right)
\cap \left( \sum_{t=0}^{\infty} W_i(s + 2t + A_{ij}) \right) \\
&\qquad\qquad\qquad\qquad\qquad\qquad\qquad\qquad \mbox{(by Lemma \ref{lem:Wsum}(i), (ii))}\\
&= \sum_{t=0}^{-A_{ij}} W_i(s + A_{ij} + 2t).
\end{align*}
as desired.
\end{proof}

\begin{lemma}
\label{lem:Z_i.Z_j}
Choose $i$ and $j$ with $1 \le i,j \le 2d_i$ and $i \neq j$. With the notation of Definition \ref{def:U.act}, the following relation holds on $V$:
\begin{equation}
\label{eq:Z_i.Z_j}
\displaystyle\sum_{r=0}^{1 - A_{ij}} (-1)^r
         \genfrac[]{0 cm}{0}{1 - A_{ij}}{r}_i 
         Z_i^{1 - A_{ij} - r} Z_j Z_i^r =
         K_i^{1- A_{ij}} K_j \prod_{s=0}^{-A_{ij}} (1 - q_i^{A_{ij} + 2s})
\end{equation}
\end{lemma}

\begin{proof}
To prove (\ref{eq:Z_i.Z_j}), it suffices to show
\begin{equation}
\label{eq:5.9}
\displaystyle\sum_{r=0}^{1 - A_{ij}} (-1)^r
         \genfrac[]{0 cm}{0}{1 - A_{ij}}{r}_i 
         Z_i^{1 - A_{ij} - r} (Z_j - K_j) Z_i^r = 0
\end{equation}
and
\begin{equation}
\label{eq:5.10}
\displaystyle\sum_{r=0}^{1 - A_{ij}} (-1)^r
         \genfrac[]{0 cm}{0}{1 - A_{ij}}{r}_i 
         Z_i^{1 - A_{ij} - r} K_j Z_i^r =
         K_i^{1- A_{ij}} K_j \prod_{s=0}^{-A_{ij}} (1 - q_i^{A_{ij} + 2s}).
\end{equation}
We first show (\ref{eq:5.9}). Since $\{W_i(s)\}_{s=0}^{2d_i}$ is a decomposition of $V$, it suffices to show that (\ref{eq:5.9}) holds on $W_i(s)$ for $0 \le s \le 2d_i$. Let $s$ be given.
By Definition \ref{def:U.act}, the left-hand side of (\ref{eq:5.9})
equals $\prod_{t=0}^{-A_{ij}} (Z_i - \varepsilon_i q_i^{s+A_{ij}+2t-d_i})(Z_j - K_j)$ on $W_i(s)$, and this product is zero in view of Lemma \ref{lem:Z_j.W_i}. We have shown (\ref{eq:5.9}). We now show (\ref{eq:5.10}).
To this end, we first recall a few identities. By \cite[p. 6]{jantzen} we find that for integers $m \ge r \ge 1$,
\begin{equation}
\label{eq:5.11}
\genfrac[]{0 cm}{0}{m}{r}_i 
+ q_i^{-m-1} \genfrac[]{0 cm}{0}{m}{r-1}_i
= q_i^{-r} \genfrac[]{0 cm}{0}{m+1}{r}_i
\qquad (1 \le i \le n).
\end{equation}
The following identity is a special case of the $q$-binomial theorem \cite[p. 236]{gasper}. For an indeterminate $\lambda$ and for an integer $m \ge 0$,
\begin{equation}
\label{eq:5.12}
\sum_{r=0}^m (-1)^{r} \genfrac[]{0 cm}{0}{m}{r}_i \lambda^{m-r}
= \prod_{r=0}^{m-1} (\lambda q_i^{m - 1 - 2r}-1) \qquad (1 \le i \le n).
\end{equation}
Using (\ref{eq:K_j}), (\ref{eq:5.11}), and induction we obtain
\begin{equation}
\label{eq:5.13}
Z_i^m K_j = \sum_{r=0}^m q_i^{(m-r)(A_{ij}+r)}
\genfrac[]{0 cm}{0}{m}{r}_i
K_i^r K_j Z_i^{m-r}
\prod_{s=0}^{r-1} (1 - q_i^{A_{ij} + 2s}).
\end{equation}
Evaluating the left-hand side of (\ref{eq:5.10}) using (\ref{eq:5.13}) we obtain
\begin{align}
\sum_{r=0}^{1-A_{ij}} &\sum_{s=0}^{1-A_{ij}-r} (-1)^r
\genfrac[]{0 cm}{0}{1-A_{ij}}{r}_i
\genfrac[]{0 cm}{0}{1-A_{ij}-r}{s}_i
q_i^{(1-A_{ij}-r-s)(A_{ij}+s)} \nonumber \\
& \times K_i^s K_j Z_i^{1 - A_{ij}}
\prod_{t=0}^{s-1} (1 - q_i^{A_{ij} + 2t}).
\label{eq:5.14}
\end{align}
For $0 \le s \le 1 - A_{ij}$ the coefficient of $K_i^s K_j Z_i^{1 - A_{ij}}$ in (\ref{eq:5.14}) is equal to
\begin{equation}
\label{eq:5.15}
\genfrac[]{0 cm}{0}{1 - A_{ij}}{s}_i
\sum_{r=0}^{1 - A_{ij} - s} (-1)^r
\genfrac[]{0 cm}{0}{1 - A_{ij} - s}{r}_i
q_i^{(1 - A_{ij} - s - r)(A_{ij}+s)}
\prod_{t=0}^{s-1} (1 - q_i^{A_{ij} + 2t}).
\end{equation}
In view of (\ref{eq:5.12}), the expression (\ref{eq:5.15}) is equal to $\delta_{s,1-A_{ij}} \prod_{t=0}^{s-1} (1 - q_i^{A_{ij} + 2t})$. Equation (\ref{eq:5.10}) follows. We have shown (\ref{eq:5.9}) and (\ref{eq:5.10}). Adding these equations we get (\ref{eq:Z_i.Z_j}) as desired.
\end{proof}

\noindent {\bf Proof of Theorem \ref{thm:Uge0}.}
To get a $\U$-module structure on $V$, we initially work with the equitable presentation for $\U$ given in Proposition \ref{prop:eqU}. Recall that in Definition \ref{def:U.act}, we defined actions of the equitable generators $K_i^{\pm 1}$, $Y_i$, and $Z_i$ $(1 \le i \le n)$ on $V$. We now show that these actions satisfy the relations of $\U$. The actions satisfy (E1)--(E3) by (e1)--(e3), (E4) by (\ref{eq:K_j}), (E5) by Lemma \ref{lem:U.rel}(ii), (E6) by (\ref{eq:Y_j}), (E7) by (e4), and (E8) by (\ref{eq:Z_i.Z_j}).  Thus these actions give a $\U$-module structure on $V$. By Definition \ref{def:U.act}, $K_i^{\pm 1} - \varepsilon_i \alpha_i^{\mp 1}k_i^{\pm 1}$ vanishes on $V$ for $1 \le i \le n$. By Remark \ref{rem:E_i.act}, $E_i - \varepsilon_i \alpha_i^{-1} e_i$ vanishes on $V$ for $1 \le i \le n$. The $\U$-module $V$ is irreducible since $V$ is an irreducible $U^{\ge 0}$-module. It is straightforward to check that $V$ has type $\varepsilon$.

Concerning the ``uniqueness'' assertion, suppose we are given a $\U$-module structure on $V$ such that the operators $E_i - \varepsilon_i \alpha_i^{-1} e_i$ and $K_i^{\pm 1} - \varepsilon_i \alpha_i^{\mp 1}k_i^{\pm 1}$ vanish on $V$. We show that this structure coincides with the one described in Definition \ref{def:U.act}. Clearly the actions of $K_i$ and $E_i$ are as in Definition \ref{def:U.act}. Thus we may describe the given structure as a set of linear operators $K_i$, $E_i$, and $\widetilde{F}_i$ $(1 \le i \le n)$ satisfying relations (R1)--(R8). In contrast with $\widetilde{F}_i$, we let
$F_i$ act on $V$ as prescribed by Definition \ref{def:U.act} and the isomorphism of Proposition \ref{prop:eqU}. We now show that $F_i = \widetilde{F}_i$ for $1 \le i \le n$.
Let $i$ be given. Consider the set $X := \{v \in V | (F_i - \widetilde{F}_i)v = 0 \}$. We argue that $X = V$. Since $V$ is an irreducible $U^{\ge 0}$-module, it suffices to prove that $X \neq 0$ and $X$ is invariant under $U^{\ge 0}$. We first show that $X$ is invariant under $U^{\ge 0}$. Note that $F_i$ and $\widetilde{F}_i$ both satisfy (R4) and (R5), so $F_i - \widetilde{F}_i$ commutes with $K_j$, $E_j$ for $1 \le j \le n$. By hypothesis, $K_j = \varepsilon_j \alpha_j^{-1} k_j$ and $E_j =  \varepsilon_j \alpha_j^{-1} e_j$ on $V$. Thus $F_i - \widetilde{F}_i$ commutes with $k_j$, $e_j$ for $1 \le j \le n$, and it follows that $X$ is invariant under $U^{\ge 0}$. We now show $X \neq 0$. Take $0 \neq v \in U_i(0)$. By (R4), $F_i v$ is in the eigenspace of $K_i$ corresponding to eigenvalue $\varepsilon_i q^{-d_i - 2}$, but this space is 0 by Lemma \ref{lem:U_i}. Similarly we have $\widetilde{F}_i v = 0$. Thus $v \in X$, and we have shown $X \neq 0$. Since $V$ is an irreducible $\U$-module, we have $X = V$ as desired. It follows that the given $\U$-module structure on $V$ is identical to the $\U$-module structure described in Definition \ref{def:U.act}. Hence, the $\U$-module described in Definition \ref{def:U.act} is unique. \hfill$\Box$

\section{The proof of Theorem \ref{thm:U}}
This section is devoted to a proof of Theorem \ref{thm:U}. We begin with a few comments about the quantum algebra $U_q(\mathfrak{sl}_2)$ and its modules.

\begin{definition} \rm \cite[p. 122]{kassel}
\label{def:sl_2}
The quantum algebra $U_q(\mathfrak{sl}_2)$ is the unital associative $\F$-algebra with generators $e$, $f$, and $k$ which satisfy the following relations:
\begin{align*}
kk^{-1} &= k^{-1} k = 1, \\
kek^{-1} &= q^2e,\\
kfk^{-1} &= q^{-2}f,\\
[e,f] &= \frac{k - k^{-1}}{q - q^{-1}}.
\end{align*}
\end{definition}

\begin{lemma} \cite[p. 128]{kassel}
\label{lem:sl_2}
Let $V$ denote a finite-dimensional irreducible $U_q(\mathfrak{sl}_2)$-module. Then there exist $\varepsilon \in \{1, -1\}$ and a basis $v_0, v_1, \ldots, v_d$ for $V$ such that $kv_i = \varepsilon q^{2i-d}v_i$ for $0 \le i \le d$, $ev_i = [i+1]v_{i+1}$ for $0 \le i \le d-1$, $ev_d = 0$, $fv_i = \varepsilon [d - i + 1]v_{i-1}$ for $1 \le i \le d$, and $fv_0 = 0$.
\end{lemma}

The proof of the next lemma is straightforward.

\begin{lemma}
\label{lem:sl_2.act}
Let $V$ denote a finite-dimensional irreducible $\U$-module. Then for $1 \le i \le n$, there exists a unique $U_q(\mathfrak{sl}_2)$-module structure on $V$ such that $e - E_i$, $f - F_i$, and $k^{\pm 1} - K_i^{\pm 1}$ vanish on $V$.
\end{lemma}

\begin{definition}
\label{def:ws} \rm
Fix a Cartan subalgebra $\mathfrak{h}$ of $\mathfrak{g}$. Let $\mathfrak{h}^*$ be the dual of $\mathfrak{h}$. Fix simple roots $\alpha_1, \ldots, \alpha_n \in \mathfrak{h}^*$ and simple coroots $\alpha_1^\vee, \ldots, \alpha_n^\vee \in \mathfrak{h}$. Then for $1 \le i,j \le n$ we have $\left<\alpha_i, \alpha_j^\vee\right> = A_{ij}$. 
Let $V$ denote a finite-dimensional irreducible $\U$-module of type $\varepsilon = \{\varepsilon_i\}_{i=1}^n$. For $\lambda \in \mathfrak{h}^*$ we define the \emph{weight space $V_\lambda$} by
\[
V_{\lambda} = \{ v \in V | K_i v = \varepsilon_i q_i^{\left<\lambda, \alpha_i^\vee\right>} \mbox{ for } 1 \le i \le n\}.
\]
We denote the \emph{set of weights of $V$} by $\Lambda = \{\lambda \in \mathfrak{h}^* | V_{\lambda} \neq 0\}$. For $1 \le i \le n$ we define linear operators $r_i: \mathfrak{h}^* \rightarrow \mathfrak{h}^*$ by
\[
r_i(\lambda) = \lambda - \left<\lambda, \alpha_i^\vee\right> \alpha_i \qquad (1 \le j \le n).
\]
Let $W$ be the group generated by $r_i$ $(1 \le i \le n)$. Then $W$ is the Weyl group associated with $\mathfrak{g}$ \cite[p. 35]{kac}.
\end{definition}

\begin{lemma} \cite[Prop. 5.1]{jantzen}
\label{lem:ws}
Let $V$ be a finite-dimensional irreducible $\U$-module with set of weights $\Lambda$. Then $V = \sum_{\lambda \in \Lambda} V_{\lambda}$ (direct sum). Moreover, $E_i$ and $F_i$ are nilpotent on $V$ $(1 \le i \le n)$.
\end{lemma}

\begin{lemma}
\label{lem:ws.act}
Let $V$ be a finite-dimensional irreducible $\U$-module with set of weights $\Lambda$.
Choose $i$ with $1 \le i \le n$, and suppose $\lambda \in \Lambda$ with $\left<\lambda, \alpha_i^\vee\right> \le 0$. Then the following hold.
\begin{enumerate}
\item The restriction of $E_i^{-\left<\lambda, \alpha_i^\vee\right>}$ to $V_{\lambda}$ is an isomorphism of vector spaces from $V_{\lambda}$ to $V_{r_i(\lambda)}$.

\item For all $v \in V_{\lambda}$, $E_i^{1 - \left<\lambda, \alpha_i^\vee\right>}v = 0$ if and only if $F_i v = 0$.

\item The restriction of $F_i^{-\left<\lambda, \alpha_i^\vee\right>}$ to $V_{r_i(\lambda)}$ is an isomorphism of vector spaces from $V_{r_i(\lambda)}$ to $V_{\lambda}$.

\item For all $v \in V_{r_i(\lambda)}$, $F_i^{1 - \left<\lambda, \alpha_i^\vee\right>} v = 0$ if and only if $E_i v = 0$.
\end{enumerate}
\end{lemma}

\begin{proof}
(i) Using Lemma \ref{lem:sl_2.act}, view $V$ as a $U_q(\mathfrak{sl}_2)$-module such that $k$ acts as $K_i$, $f$ acts as $F_i$, and $e$ acts as $E_i$. As a $U_q(\mathfrak{sl}_2)$-module, $V$ is the direct sum of irreducible $U_q(\mathfrak{sl}_2)$-modules. Let $S$ be an irreducible $U_q(\mathfrak{sl}_2)$-module summand with $S \cap V_\lambda \neq 0$. Then by Lemma \ref{lem:sl_2}, $E_i^{-\left<\lambda, \alpha_i^\vee\right>} = e^{-\left<\lambda, \alpha_i^\vee\right>}$ is an isomorphism from $S \cap V_\lambda$ to $S \cap V_{r_i(\lambda)}$. The result follows.

(ii) Without loss of generality, $v \in S$ for some irreducible $U_q(\mathfrak{sl}_2)$-submodule $S$ of $V$. Then $E_i^{1 - \left<\lambda, \alpha_i^\vee\right>} v = e^{1 -\left<\lambda, \alpha_i^\vee\right>} v = 0$ if and only if $F_i v = f v = 0$ by Lemma \ref{lem:sl_2}. The result follows.

(iii), (iv) These results follow in a manner analogous to (i), (ii).
\end{proof}

\begin{corollary}
\label{cor:Weyl}
Let $V$ be a finite-dimensional irreducible $\U$-module with set of weights $\Lambda$. Let $W$ be the Weyl group from Definition \ref{def:ws}. Then for $w \in W$, $V_\lambda \neq 0$ if and only if $V_{w(\lambda)} \neq 0$.
\end{corollary}

\begin{proof}
Since $W$ is generated by the linear operators $r_i$, we may assume $w = r_i$ for some $1 \le i \le n$. The result now follows from Lemma \ref{lem:ws.act}(i), (iii).
\end{proof}

\begin{lemma}
\label{lem:Lambda}
Let $V$ be a finite-dimensional irreducible $\U$-module with set of weights $\Lambda$. Choose $\lambda \in \Lambda$ with $\lambda \neq 0$. Then there exist integers $i$ and $j$ with $1 \le i, j \le n$ such that $\left<\lambda, \alpha_i^\vee \right> < 0$ and $\left<\lambda, \alpha_j^\vee \right> > 0$.
\end{lemma}

\begin{proof} Recall the integers $u_i$ from Definition \ref{def:A}. Define $K \in \mathfrak{h}$ by $K := \sum_{i=1}^n u_i \alpha_i^\vee$ (cf. \cite[p. 80]{kac}). Since $\lambda \neq 0$ and $u_i > 0$ for $1 \le i \le n$, it suffices to show $\left< \lambda, K \right> = 0$. To this end, we define $C \in \U$ by $C := \prod_{i=1}^n K_i^{u_i}$. Note that for any $v \in V_{\lambda}$, we have $Cv = \pm q_i^{\left<\lambda, K\right>} v$. We claim that $C$ acts on $V$ as $\pm 1$. It will follow that $\left<\lambda, K\right> = 0$ as desired. 
We observe that $C$ is the central element of $V$ described in \cite[Thm. 12.2.1]{chari}. By \cite[Prop. 12.2.3]{chari}, $C$ acts on $V$ as $\pm 1$, so the claim holds. The result follows.
\end{proof}

\begin{lemma}
\label{lem:W}
Let $V$ be a finite-dimensional irreducible $\U$-module with set of weights $\Lambda$. Let $S$ be a subspace of $V$ such that $K_i S \subseteq S$ and $E_i S \subseteq S$ for $1 \le i \le n$. Then we have the following:
\begin{enumerate}
\item $S = \sum_{\lambda \in \Lambda} (S \cap V_{\lambda})$,

\item For $1 \le j \le n$, $S \cap V_{\lambda} \neq 0$ if and only if $S \cap V_{r_j(\lambda)} \neq 0$.
\end{enumerate}
\end{lemma}

\begin{proof}
Since $S$ is invariant under $K_i$ $(1 \le i \le n)$, we have (i). To see (ii), note that if $\left<\lambda, \alpha_j^\vee \right> \le 0$ then $E_j^{-\left<\lambda, \alpha_j^\vee \right>}$ is an isomorphism from $S \cap V_{\lambda}$ to $S \cap V_{r_j(\lambda)}$ by Lemma \ref{lem:ws.act}(i) and the hypothesis $E_j S \subseteq S$, so the result holds in this case. If $\left<\lambda, \alpha_j^\vee \right> > 0$ we argue as follows.

Let $W\lambda$ be the orbit of $\lambda$ under $W$. $V$ is finite-dimensional, so $W\lambda$ is a finite set by Lemma \ref{cor:Weyl}. Construct a directed graph with vertex set $W\lambda$ and edge set consisting of all edges $(\mu_1, \mu_2) \in W\lambda \times W\lambda$ with the following property:
\begin{equation}
\label{eq:graph}
\mbox{there exists $h$ $(1 \le h \le n)$ such that $\left<\mu_1, \alpha_h^\vee \right> < 0$ and $\mu_2 = r_h(\mu_1)$}.
\end{equation}
If $(\mu_1, \mu_2)$ is an edge of the graph, (\ref{eq:graph}) guarantees $S \cap V_{\mu_1} \cong S \cap V_{\mu_2}$ in view of Lemma \ref{lem:ws.act}(i) and the hypothesis $E_i S \subseteq S$ $(1 \le i \le n)$. Moreover if $\mu_1, \mu_2 \in W\lambda$ and there exists a path from $\mu_1$ to $\mu_2$ in the graph, we have $S \cap V_{\mu_1} \cong S \cap V_{\mu_2}$. To complete the proof, we show that there is a path from $\lambda$ to $r_j(\lambda)$.

Note that $0 \notin W\lambda$, since otherwise $W\lambda$ contains only the element $0$ in contradiction of $\left<\lambda, \alpha_j^\vee \right> > 0$. Thus for each $\mu \in W\lambda$, there are integers $h$ and $h'$ with $1 \le h, h' \le n$ such that $\left< \mu, \alpha_h^\vee \right> < 0$ and $\left< \mu, \alpha_{h'}^\vee \right> > 0$ by Lemma \ref{lem:Lambda}. Therefore, by the construction of the graph, each vertex has at least one incoming and at least one outgoing edge, and no vertex is self-adjacent. 
At each vertex, color exactly one incoming and exactly one outgoing edge.
We have $\left<r_j(\lambda), \alpha_j^\vee \right> < 0$, so we may assume the edge $(r_j(\lambda), \lambda)$ is colored.
A simple argument shows that $\lambda$ is contained in a cycle of colored edges. Thus there is a path from $\lambda$ to $r_j(\lambda)$. The result follows.
\end{proof}

\noindent {\bf Proof of Theorem \ref{thm:U}.} Let $V$ be a finite-dimensional irreducible $\U$-module of type $\varepsilon$. Suppose $\alpha \in (\F^\times)^n$. We first prove that the desired $U^{\ge 0}$-module structure on $V$ exists. For $1 \le i \le n$ let $e_i$ and $k_i^{\pm 1}$ act on $V$ as $\varepsilon_i \alpha_i^{-1} E_i$ and $\varepsilon_i \alpha_i^{\mp 1} K_i^{\pm 1}$ respectively. Then using the defining relations for $\U$ in Definition \ref{def:U}, it is easy to see that $e_i$ and $k_i^{\pm 1}$ satisfy (r1)--(r4) and therefore induce a $U^{\ge 0}$-module structure on $V$. From the construction, the operators $E_i - \varepsilon_i \alpha_i^{-1} e_i$ and $K_i^{\pm 1} - \varepsilon_i \alpha_i^{\mp 1}k_i$ vanish on $V$. We have now shown the desired $U^{\ge 0}$-module structure exists, and it is clear that this $U^{\ge 0}$-module structure is unique. Next we show that the $U^{\ge 0}$-module structure is irreducible. To this end, we let $S$ denote an irreducible $U^{\ge 0}$-submodule of $V$ and argue that $S = V$. Since $V$ is an irreducible $\U$-module, it suffices to show that $S$ is nonzero and invariant under $\U$. By construction, $S$ is nonzero and invariant under $E_i, K_i^{\pm 1}$ $(1 \le i \le n)$. We claim that $S$ is invariant under $F_i$ $(1 \le i \le n)$. 
Let $i$ be given, and define $\widetilde{S} := \{v \in S | F_i v \in S \}$. To prove the claim, it suffices to show $S = \widetilde{S}$. Since $S$ is an irreducible $U^{\ge 0}$-module, it is enough to show that $\widetilde{S}$ is nonzero and invariant under $U^{\ge 0}$. Fix $j$ with $1 \le j \le n$. It follows from (R5) and the fact that $E_j - \varepsilon_j \alpha_j^{-1} e_j$ and $K_j^{\pm 1} - \varepsilon_j \alpha_j^{\mp 1}k_j$ are $0$ on $V$ that $e_j \widetilde{S} \subseteq \widetilde{S}$.  By (R3) we have $K_j E_i = q_j^{A_{ji}} E_i K_j$. By this and since $K_j - \varepsilon_j \alpha_j^{-1}k_j$ is $0$ on $V$, we find that $k_j E_i - q_j^{A_{ji}} E_i k_j$ vanishes on $V$. Consequently, we have $k_j \widetilde{S} \subseteq \widetilde{S}$, and then $k_j^{-1} \widetilde{S} \subseteq \widetilde{S}$ holds as well. Thus $\widetilde{S}$ is invariant under $U^{\ge 0}$.
To verify that $\widetilde{S} \neq 0$, observe that by Lemma \ref{lem:W}(i), we have $S = \sum_{\lambda \in \Lambda} (S \cap V_{\lambda})$. By Lemma \ref{lem:ws}, $E_i$ is nilpotent on $V$, so we can choose $\lambda \in \Lambda$ such that $S \cap V_{\lambda} \neq 0$ and $E_i (S \cap V_{\lambda}) = 0$. In view of Lemma \ref{lem:ws.act}(i) and our choice of $\lambda$, we have $\left<r_i(\lambda), \alpha_i^\vee \right> \le 0$. By Lemma \ref{lem:W}(ii), we can choose $0 \neq v \in S \cap V_{r_i(\lambda)}$. It follows that $E_i^{ 1 - \left< r_i(\lambda), \alpha_i^\vee \right> } v \in E_i (S \cap V_\lambda) = 0$. Thus by Lemma \ref{lem:ws.act}(ii), $F_i v = 0$. We have $v \in \widetilde{S}$, so $\widetilde{S} \neq 0$. We have shown that $\widetilde{S}$ is nonzero and invariant under $U^{\ge 0}$. Thus $\widetilde{S}=S$ since $S$ is an irreducible $U^{\ge 0}$-module. The claim follows. We have shown that $S$ is nonzero and invariant under $\U$, so $S=V$ since $V$ is an irreducible $\U$-module. Thus $V$ is an irreducible $U^{\ge 0}$-module as desired. It is routine to show the $U^{\ge 0}$-module structure on $V$ has type $\alpha$.
\hfill$\Box$

\section{Irreducible $\U^{\ge 0}$-modules}

\begin{theorem}
\label{thm:borel}
Suppose $\varepsilon \in \{1, -1\}^n$. Then the following hold:
\begin{enumerate}
\item Let $V$ be a finite-dimensional irreducible $\U^{\ge 0}$-module of type $\varepsilon$. Then the action of $\U^{\ge 0}$ on $V$ extends uniquely to an action of $\U$ on $V$.
The resulting $\U$-module structure on $V$ is irreducible and of type $\varepsilon$.

\item Let $V$ be a finite-dimensional irreducible $\U$-module of type $\varepsilon$. 
When the $\U$-action is restricted to $\U^{\ge 0}$, the resulting $\U^{\ge 0}$-module
structure on $V$ is irreducible and of type $\varepsilon$.
\end{enumerate}
\end{theorem}

\begin{proof}
(i) We define an action of $U^{\ge 0}$ on $V$ such that $k_i^{\pm 1}$ acts as $K_i^{\pm 1}$ and $e_i$ acts as $E_i$ for $1 \le i \le n$. This $U^{\ge 0}$-module structure has type $\varepsilon$. By Theorem \ref{thm:Uge0}, there exists a unique $\U$-module structure on $V$ such that the operators $k_i^{\pm 1} - K_i^{\pm 1}$ and $E_i - e_i$ vanish on $V$ for $1 \le i \le n$. This structure is irreducible and of type $\varepsilon$, and it restricts to the $\U^{\ge 0}$-module structure on $V$. The result follows.

(ii) We let $S$ denote a nonzero $\U^{\ge 0}$-submodule of $V$ and claim that $S = V$. By its definition, $S$ is invariant under $K_i^{\pm 1}, E_i$ for $1 \le i \le n$. By Theorem \ref{thm:U}, for any $\alpha \in (\F^\times)^n$, there is the structure of a $U^{\ge 0}$-module on $V$ such that the operators $K_i^{\pm 1} - \varepsilon_i \alpha_i^{\mp 1}$ and $E_i - \varepsilon_i \alpha_i^{-1} e_i$ vanish on $V$ $(1 \le i \le n)$. From this we see that $S$ is invariant under $U^{\ge 0}$. But the $U^{\ge 0}$-module structure on $V$ is irreducible by Theorem \ref{thm:U}, so $S = V$ and our claim is proved. Note that $V$ has type $\varepsilon$ as a $\U^{\ge 0}$-module since it has type $\varepsilon$ as a $\U$-module.
\end{proof}

\section{Acknowledgements}
The author would like to thank Paul Terwilliger and Georgia Benkart for their helpful comments and suggestions.

\normalsize
\noindent \textsc{Department of Mathematics, University of Wisconsin, Madison, WI 53706-1388 USA} \\
{\tt bowman@math.wisc.edu}

\end{document}